%% file: cbp.tex
%

\input def1200.tex

\input mssymb

\null
\overfullrule=0pt

\def\vta{{\vtop{\offinterlineskip{\hbox{$A$}
\hbox{$\scriptstyle\sim$}}}}}
\def\vtb{{\vtop{\offinterlineskip{\hbox{$B$}\hbox
{$\scriptstyle \sim$}}}}}
\def\vtd{{\vtop{\offinterlineskip{\hbox{$D$}\hbox
{$\scriptstyle \sim$}}}}}
\def\vtj{{\vtop{\offinterlineskip{\hbox{$j$}\hbox
{$\scriptstyle\sim$}}}}}
\def\vtf{{\vtop{\offinterlineskip{\hbox{$f$}\hbox
{$\scriptstyle\sim$}}}}}

\vskip2truecm
{\nopagenumbers
\overfullrule=0pt
\sect{ON DENSITIES OF BOX PRODUCTS}
\bigskip
\ce{by}
$$\vbox{\halign{\tabskip2em\hfil#\hfil&\hfil#\hfil\cr
\bf Moti Gitik&\bf Saharon Shelah\cr
School of Mathematical Sciences&Hebrew University of
Jerusalem\cr
Sackler Faculty of Exact Sciences&Department of
Mathematics\cr
Tel Aviv University&Givat Ram, Jerusalem\cr
Ramat Aviv  69978 Israel&Israel\cr}}$$

\vskip1truecm
\dspace
{\narrower\medskip
\subheading{Abstract}

We construct two universes $V_1, V_2$ satisfying the
following GCH below $\aleph_\ome$, 
$2^{\aleph_\ome}=\aleph_{\ome +2}$ and the topological
density of the space
${}^{\aleph_\ome}\,2$ with $\aleph_0$ box
product topology $d_{<\aleph_1}(\aleph_\ome)$ is
$\aleph_{\ome+1}$  in $V_1$  and $\aleph_{\ome+2}$
in $V_2$.  Further related results are discussed
as well.\medskip}
\vfill\eject}
\count0=1
\dspace
\def\llvdash{\mathop{\|\hskip-2pt
\raise 3pt\hbox{\vrule height 0.25pt width 0.50cm}}}
\def\lvdash{\mathop{|\hskip-2pt \raise 3pt\hbox{\vrule
height 0.25pt width 0.50cm}}}
\def\l{{\langle}}
\def\smallfrown{{\scriptscriptstyle \cap}}
\def\r{{\rangle}}
\def\Ult{{\rm Ult}}
\def\ult{{\rm ult}}

\dspace

W.~Comfort asked the following question:
Assume $\lam$  is a strong limit singular,
$\kap>\cf\lam$.  Is $d_{<\kap}(\lam)=2^\lam$? Is it
always $>\lam^+$ when $2^\lam >\lam^+$?  

$d_{<\kap}(\lam)$ denotes the density of the topological
space ${}^\lam\!2$ with topology generated by the
following family of clopen sets: 
$$\{ [f]\mid f\in {}^{\! a} 2\quad \hbox{for some}\quad
a\subseteq\lam, |a|<\kap\}$$
where $[f]=\{ g\in{}^{\!\lam} 2\mid g\supseteq f\}$,
i.e. $d_{<\kap}(\lam)=\min\{|F|\mid F\subseteq
{}^{\!\lam}2$  and if $a\subseteq \lam$
$|a|<\kap$  and $g\in{}^{\! a}2$  then there is
$f\in F$  $g\subseteq f\}$.  

The aim of this paper will be to show that under
$\neg SCH$ $d_{<\aleph_1}(\lam)$  may be
$\lam^+$  even if $2^\lam >\lam^+$.
Surprisingly, it turned out that it is easier to
get $d_{<\aleph_1}(\lam)=\lam^+$  than
$d_{<\aleph_1}(\lam)=2^\lam$  for a strong limit
$\lam$  of cofinality $\aleph_0$  with $2^\lam >\lam^+$.
We refer to the ZFC results using the cardinal
arithmetic to Shelah [Sh430, \S5].

The paper is organized as follows. Section 1
is less involved and provides a model with a
strong limit $\lam$,  $cf\lam =\aleph_0$,
$2^\lam >\lam^+$ and $d_{<\aleph_1}(\lam)$.  The
main disadvantage is that $\lam$  is rather
large and it is unclear how to move everything 
down to say $\aleph_\ome$.  But as a bonus this
construction gives a normal ultrafilter over
$\lam$  generated by $\lam^+$  sets and $2^\lam
>\lam^+$.  Originally such models were produced
by T.~Carlson and H.~Woodin (both unpublished).
In Section 2 it is fixed by the cost of using more
involved techniques.  Also initial assumptions reduced
from huge to hypermeasurable.

Both section can be read independently.  Most of the
construction in Section 1 is due to the second author.
Only the final argument using a huge cardinal
is of the first author.  The construction in
Section 2 is due to the first author.  

\sect{1.~~Density of Box Products From
Huge Cardinal}

In this section, we prove the following:

\proclaim Theorem 1.1.  Suppose that $\lam$  is
a huge cardinal. Then there exists a
generic extension satisfying the following:
\item{(a)} $\lam$  is a strong limit of cofinality
$\ome$.
\smallskip
\item{(b)} $2^\lam >\lam^+$.
\smallskip
\item{(c)} for every $\mu <\lam$ $d_{<\mu}(\lam)=\lam^+$.

Let $\kap$  be a measurable cardinal.  Let $D$
be a normal ultrafilter over $\kap$.

\subheading{Definition 1.2}  Let $Q_D$  be a
forcing notion consisting of all triples $\langle
f,\alp,A\rangle$ so that 
\item{(a)} $A\in D$
\smallskip
\item{(b)} $\alp <\kap$
\smallskip
\item{(c)} $f$  is a function on $[A]^{<\ome}$
such that
\smallskip
\item{(c1)} for every $\eta\in [A]^{<\ome}$
$f(\eta)$  is a partial function from min
$(A\bks (\max\eta +1))$  to $2$
\smallskip
\item{(c2)} $\sup \{ |\dom f(\eta )|\mid \eta\in
[A]^{<\ome}\} <\kap$.

\subheading{Definition 1.3} Let $\langle
f_1,\alp_1,A_1\rangle$, $\langle f_2,\alp_2,A_2\rangle\in
Q_D$.  We define $\langle f_2,\alp_2,A_2\rangle
\ge\langle f_1,\alp_1,A_1\rangle$ iff
\item{(a)} $\alp_1\le\alp_2$
\smallskip
\item{(b)} $A_1\supseteq A_2$
\smallskip
\item{(c)} $A_1\cap\alp_1=A_2\cap\alp_1$
\smallskip
\item{(d)} for every $\eta\in [A_2]^{<\ome}$
$f_1(\eta)\subseteq f_2(\eta)$
\smallskip
\item{(e)} $f_1\rhookup [A_1\cap\alp_1]^{<\ome}=f_2
\rhookup [A_2\cap\alp_1]^{<\ome}$.
\item{}

\smallskip
Intuitively, the forcing is intended to add a
set $A\subseteq\kap$  which is almost contained
in every set of $D$  and a function $f$ on
$[A]^{<\ome}$  which is a name of a function in
a Prikry forcing for changing cofinality of
$\kap$  to $\aleph_0$.  This function will be
eventually a member of a desired dense set of
cardinality $\kap^+$.

The idea will be to add $\lam$  new subsets to
$\kap$  $(\lam=\kap^{++}$  or any desired value
for the final $2^\kap$) preserving
supercompactness of $\kap$  together with
iteration of the length $\kap^+$  of forcings
$Q_{\vtop{\offinterlineskip\hbox{$D_i$}\hbox{$\sim$}}}
(i<\kap^+)$,  where 
${\vtop{\offinterlineskip\hbox{$D_i$}\hbox{$\sim$}}}$'s
are picked to increase. Finally we'll obtain $\buildrul
\sim\under D=\cup {\vtop{\offinterlineskip\hbox{$D_i$}
\hbox{$\sim$}}}$ and force with the Prikry
forcing for $\vtd$.  The
interpretation of the generic functions $f_i$'s
$(i<\kap^+)$ from each stage of the iteration
will form the dense set of cardinality
$\kap^+$.

Let us start with a basic fact about names in
the Prikry forcing.

\proclaim Lemma 1.4.  Let $D$  be a normal ultrafilter
over $\kap$, $\calP_D$ the Prikry forcing with
$D$, $\tau$  a $\calP_D$-name
of a partial function of cardinality $<\mu (\mu <\kap)$
from $\kap$  to $2$.  Then there are $A$  and
$f$  satisfying the conditions (a), (c) of 1.1 so that 
$\langle\phi, A\rangle\llvdash\buildrul\sim\under\tau
= \bigcup_{n<\ome}f((
{\vtop{\offinterlineskip\hbox{$\kap$}\hbox{$\scriptstyle
\sim$}}}{}_n, {\vtop{\offinterlineskip\hbox{$\kap$}
\hbox{$\sim$}}}{}_1, \nek{\vtop{\offinterlineskip\hbox
{$\kap$} \hbox{$\scriptstyle\sim$}}}{}_n))$ where
$\langle{\vtop{\offinterlineskip\hbox{$\kap$}
\hbox{$\scriptstyle\sim$}}}{}_n\mid
n<\ome\rangle$  is the canonical name of the
Prikry sequence.  Also $|f(\eta)|<\mu$ for each
$\eta \in [A]^{<\ome}$.  

\pr Using normality, we pick $A\in D$  and
$\langle | a_\eta|\eta\in [A]^{<\ome}\rangle$,
$|a_\eta |<\mu$  $(\eta\in [A]^{<\ome})$
such that for every $\eta\in [A]^{<\ome}$  $\langle
\eta, A\bks \max\eta\rangle$  $\llvdash\buildrul\sim\under
\tau\cap (\max\eta$, the first element of the
Prikry sequence above $\eta$) $=a_\eta$.

Define $f(\eta)=a_\eta$  for $\eta\in
[A]^{<\ome}$.  Then, clearly
$$\langle\phi,A>\llvdash\buildrul\sim\under\tau
=\bigcup_{n<\ome}f({\vtop{\offinterlineskip\hbox{$\kap$}
\hbox{$\scriptstyle\sim$}}}{}_0\nek{\vtop{
\offinterlineskip\hbox{$\kap$} \hbox{$\scriptstyle\sim$}}}
{}_n))$$
\hfill$\square$

Let $G\subseteq Q_D$.  We define $A_D=\cap \{
A\in D|$ for some $\alp,f$ $\langle f,\alp,
A\rangle \in G\}$  and $f_{D,\mu}$  will be a
function with domain $[A_D]^{<\ome}$  so that
for every $\eta\in [A_D]^{<\ome}$  $f_D(\eta)
=\cup \{ f(\eta)|$  for some $\alp,A$ $\langle
f,\alp, A\rangle\in G\}$.  

Let ${\vtop{\offinterlineskip\hbox{$A$}\hbox
{$\scriptstyle\sim$}}}_D$, ${\vtop{\offinterlineskip
\hbox{$f$}\hbox{$\scriptstyle\sim$}}}_D$, be a
canonical name of $A_D$, $f_D$.

Let $\calP_D$  denote the Prikry forcing with $D$.

The following lemma is crucial.

\proclaim Lemma 1.5.  Suppose that $D$ is a
normal ultrafilter over $\kap$  and $\buildrul\sim
\under\tau$  is a $\calP_D$-name of partial
function of cardinality $<\mu$  (for some $\mu
<\kap$) from $\kap$ to $2$.

{\sl Suppose that $\langle\phi,0,\kap\rangle
\llvdash_{Q_D}$ ``there is a normal ultrafilter
${\vtop{\offinterlineskip\hbox{$D$}\hbox
{$\scriptstyle\sim$}}}_1$, over $\kap$  with
${\vtop{\offinterlineskip\hbox{$A$}\hbox{$\scriptstyle
\sim$}}}_D\in {\vtop{\offinterlineskip\hbox{$D$}\hbox
{$\scriptstyle\sim$}}}_1$".
Then there is a generic $G\subseteq Q_D$, so
that if $D_1$  is a normal ultrafilter in $V[G]$
with $A_G\in D_1$, {\it then\/}, in $V[G]$}
$$\langle \phi, A_G\rangle\llvdash_{\calP_{D_1}}
\buildrul\sim\under\tau\subseteq\bigcup_{n<\ome}f_G
(\langle \kap_0\nek \kap_n\rangle )\ .$$

\pr Applying Lemma 1.4 to $D,\tau$  in $V$  we
pick $A$,  $f$  as in the conclusion of the
lemma.  Now let $G\subseteq Q_D$ be generic with
$\langle f,0,A\rangle\in G$.

Then $A_G\subseteq A$  and for every
$\eta\in [A_G]^{<\ome}$  $f(\eta)\subseteq f_G(\eta)$,
by 1.2.  But since  
$$\langle \phi,A\rangle \llvdash_{\calP_D}\buildrul
\sim\under\tau =\bigcup_{n<\ome}f(\langle
{\vtop{\offinterlineskip\hbox{$\kap$}
\hbox{$\scriptstyle \sim$}}}_0\nek
{\vtop{\offinterlineskip\hbox{$\kap$}
\hbox{$\sim$}}}{}_n\rangle)$$
and $D\subseteq D_1$  we are done.\hfill$\square$

Now the plan will be as follows:
We'll blow up the power of $\kap$  to some
cardinal of cofinality $\kap^+$ using
$<\kap$-support iteration of forcings of the
type $Q_\vtd$.  Using hugeness, a sequence
$$D_0\subseteq\vtd_1\subseteq\vtd_2\subseteq
\cdots\subseteq\vtd_\alp\subseteq\cdots
(\alp <\kap^+)$$
will be generated and $Q_{\vtd_{\alp}}$'s will
be used cofinally.  The final step will be to
use the Prikry forcing with $\bigcup_{\alp
<\kap^+}\vtd_{{}\alp}$.

Let us observe first that the forcing $Q_D$ is quite nice.

\proclaim Lemma 1.6.  $Q_D$  is $<\kap$-directed complete.

\pr It is obvious from the definition.

\subheading{Definition 1.7} (Shelah [Sh80]).
Let $P$  be a forcing notion.  $P$  satisfies a
``stationary" $\kap^+$-c.c. iff for every
$\langle p_i\mid i<\kap^+\rangle$  in $\calP$
there is a closed unbounded set $C\subseteq\kap^+$
and a regressive function $f:\kap^+\to\kap^+$
such that for $\alp,\bet\in C$  if $cf\alp =cf
\bet =\kap$  and $f(\alp)=f(\bet)$  then
$p_\alp$  and $p_\bet$  are compatible.

\proclaim Lemma 1.8. $Q_D$  satisfies the
``stationary" $\kap^+$-c.c.

\pr Let $p_i=\langle f_i,\alp_i, A_i\rangle$
$(i<\kap^+)$  be conditions in $Q_D$.

For every $\sig,\alp <\kap$, $a\subseteq
\alp$, $g$  a function on $[a]^{<\ome}$  we
set $A_{\sig, \alp, a,g}=\{ i<\kap^+\mid\sig =\sup
\{ |\dom f_i(\eta)|\eta\in [A_i]^{<\ome}\}$,
$\alp=\alp_i$, $A_i\cap\alp_i=a$  and
$f_i\rhookup [a]^{<\ome}=g\}$.

Then $\kap^+$  is a disjoint union of these
$\kap$  sets.

It is enough to prove the following claim:

\subheading{Claim} For every $\sig,\alp,
a,g$  as above among any $\big(2^{|\sig|+|\alp|}
\big)^+$ members of $A_{\sig,\alp,a,g}$  at
least two are compatible.

Let us first complete the proof of the lemma
using the claim and then we prove the claim.

Denote $A_{\sig,\alp, a,g}$  by $A$.  Assume
that $\{\del <\kap^+|cf\del =\kap\}\cap A$
is stationary.  Clearly, there are $\sig,\alp,
a,g$ for which this is true.  Let $\del\in
A$,  $cf\del =\kap$.  We define by induction
on $\eps$ an increasing sequence of ordinals
$\alp_{\del,\eps}<\del$  in $A$  such that
$p_{\alp_{\del,\eps}}$  is incompatible with
$p_\del$  and with $p_{\alp_{\del,\rho}}$
for $\rho <\eps$.  At stage $\eps$  just
pick $\alp <\del,\alp\in A$  such that
$p_\alp$ is incompatible with $p_\del$  and
every $p_{\alp_{\del,\rho}}(\rho <\eps)$ if
there is such an $\alp$.  Otherwise we stop.
Let $\langle \alp_{\del,\eps}\mid\eps
<\tau_\del\rangle$  be such a sequence.
Then, by the claim, $\tau_\del <(2^{|\alp|
+|\sig|})^+ <\kap$.  Hence, if we take a
regressive function $g(\del)=\tau_\del$,
then whenever $g(\del_1)=g(\del_2)$
$p_{\del_1},p_{\del_2}$  will be compatible.
So, we obtain a ``stationary" $\kap^+$-c.c.

\subheading{Proof of the Claim}
Let $\langle i_\xi\mid\xi <(2^{|\sig|+|\alp |})^+
\rangle$  be a sequence from $A$.  Set $B_0=
\bigcap_{\xi <2^{|\sig |+|\alp|}} A_{i_\xi}$. 
Then $B_0\in D$.  There is
$B_1\subseteq B_0$,  $B_1\in D$  such that
the isomorphism types of structures
$$\langle\alp,\rho ,a,\langle\dom f_{i_\xi}(\nu^\cap\rho
\rhookup\ell)|\xi <(2^{|\sig|+|\alp|})^+,\nu\in
[a]^{<\ome},\ell\le\ {\rm length}(\rho
\rangle , \le\rangle$$
depends only on the length
of $\rho$  for $\rho\in [B_1]^{<\ome}$.  
Choose $\eps_0<\eps_1<\cdots
<\eps_n<\cdots (n<\ome)$  an $\ome$-sequence
of elements of $B_1$.  Now using Erd\"os-Rado
Theorem it is easy to find $\xi_0
<\xi_1<(2^{|\sig|+|\alp|})^+$  such that for
every $\rho\in [a\cup \{\eps_\ell|\ell
<\ome\}]^{<\ome}$  $f_{i_{\xi_0}}(\rho)$
and $f_{i_{\xi_1}}(\rho)$  are compatible.
But then $f_{i_{\xi_0}}(\rho)$  and $f_{i_{\xi_1}}
(\rho)$  will be compatible for every
$\rho\in [B_1]^{<\ome}$. Which implies a
compatibility of $p_{i_{\xi_0}}$ and
$p_{i_{\xi_1}}$.\hfill $\square$ of the
claim.

Let now $\kap$  be an almost huge cardinal
with a measurable target point, i.e. there is
$j:V\to M$,  critical $(j)=\kap$,  $j(\kap)=\lam$,
${}^{\lam >}\! M\subseteq M$  and $\lam$
is a measurable cardinal in $V$.  Fix such
an embedding $j:V\to M$,  $j(\kap)=\lam$
and a normal measure $U_\lam$  over $\lam$.

We define an iteration
$\langle P_\alp,Q_{0\alp}*Q_{1\alp}\mid\alp
<\kap\rangle$  as follows:  if $\alp$  is
not measurable in $V^{P_\alp}$  then
$Q_{0\alp}*Q_{1\alp} =\emptyset$; if $\alp$
is a measurable cardinal in $V^{P_\alp}$,
then $Q_{0\alp}$ will be atomic forcing
picking an ordinal $F(\alp)<\kap$  and
$Q_{1\alp}$  will $<\alp$-support iteration
of the maximal possible length $\le F(\alp)$
of forcings of the form
$Q_{\vtop{\offinterlineskip\hbox{$D$}
\hbox{$\scriptstyle\sim$}}}$ over all
normal ultrafilters
${\vtop{\offinterlineskip\hbox{$D$}
\hbox{$\scriptstyle\sim$}}}$ over
$\alp$.  I.e. first over $\alp$  we force with
$<\alp$-supported product of forcings $Q_D$ where
$\buildrul\sim\under D$  runs over all normal
ultrafilters over $\alp$.  If $\alp$  remains
measurable after this forcing, then again we
force with $Q_D$'s for each normal ultrafilter $D$ 
of this extension and so on as as far as possible up to
$F(\alp)$.  Easton support is used at limit stages
of the iteration.  By Shelah [Sh80] and Lemmas
1.8, 1.6 $Q_{1\alp}$  satisfies $\alp^+$-c.c.
and is $\alp$-directed closed over $V^{\calP_\alp*
Q_{0\alp}}$, for $\alp <\kap$. 

The role of the trivial forcing $Q_{0\alp}$  is
to bound the length of the iteration of
$Q_{1\alp}$.  It is needed, since, for example,
if $\alp$  is a supercompact and $\alp$-directed
closed indestructible, then forcings $Q_D$ will
preserve its supercompactness and hence also the
measurability.  So new ultrafilters will appear
over $\alp$  forever.  

Let us work now over $\kap$.  Let
$G_\kap\subseteq P_\kap$  be generic.  We
consider in $M$  $j(P_\kap)=P_{j(\kap)}$  and
$P_{j(\kap)}\big/ G_\kap$  in $M[G_\kap]$.  Let
us split $P_{j(\kap)}\big/ G_\kap$  into $Q_{0\kap}*
Q_{1\kap}$  and $P_{>\kap}$.  The generic object
for $Q_{0\kap}$ is just any ordinal
$F(\kap)<j(\kap)=\lam$.  By standard arguments
on backwards Easton forcing (see, for example, A.
Kamamori [Ka]), for every $F(\kap)\le\lam$  the
length of $Q_{1\kap}$  will be $F(\kap )$. 
For a while set  $F(\kap)=\lam$,  i.e. we like
to deal with iteration $Q_{1\kap}$ of the length
$\lam$.  We consider an enumeration $\langle
{\vtop{\offinterlineskip\hbox{$A$}
\hbox{$\scriptstyle\sim$}}}{}_\tau\mid\tau <\lam
\rangle$  of  $Q_{1\kap}$-names of subsets of
$\kap$  in $M[G_\kap] [\{\lam\}]$,  such that
$\tau_1 <\tau_2<\lam$  implies that
${\vtop{\offinterlineskip\hbox{$A$}\hbox{$\scriptstyle
\sim$}}}{}_{\tau_1}$. 
depends on the part of $Q_{1\kap}$
of the length $\le$  then those of
${\vtop{\offinterlineskip\hbox{$A$}\hbox{$\scriptstyle
\sim$}}}{}_{\tau_2}$. 
Since $\lam$  is measurable and
$Q_{1\kap}$  has $<\kap$-support there will be
$C\in \calU_\lam$ consisting of inaccessibles such that
for every
$\del\in C\langle\vtop{\offinterlineskip\hbox{$A$}
\hbox{$\scriptstyle\sim$}}{}_\tau\mid\tau <\del\rangle$
enumerates the names of all subsets of $\kap$
appearing before the stage $\del$, i.e.
$Q_{1\kap}\rhookup\del$-names.  Equivalently, all the
subsets for the $Q_{1\kap}$ with $F(\kap)=\del$.  Now let
$\del$  be in $C$.  For every $\buildrul\sim\under D$
appearing in $Q_{1\kap}\rhookup \del$  let
$r_{{}\!\buildrul\sim\under D}\in Q_{j(\buildrul\sim\under
D)}$ be defined as follows.  $r_{{}\!\buildrul\sim\under
D}= \langle\buildrul\sim\under f,\kap ,\buildrul\sim\under
A\rangle$  where 
$$\vta=A_{\vtd}\cup(\cap \{ \vtj
(\vtb)\mid\vtb\in\vtd\})\ ,$$
$\vtf \rhookup \kap =f_{\vtd}$  and above $\kap$
we take $\vtf(\eta)=\bigcup \{ j(\vtf
)(\eta)\mid \vtf$  appear in a condition in $G_\vtd\}$.

Let $q_\del\in Q_{1j(\kap)}$  consists of this
$r_\vtd$'s sitting in the right place.

Clearly, that if $\rho >\del$  is also in $C$,
then $q_\rho\rhookup\del=q_\del$.   

Let $\rho$  be in $C$.  Pick a master condition $p_{\rho}
\in P_{>\kap}*(j(\rho)*Q_{1j(\kap)})$  deciding all
the statements ``$\kap\in j(\vta_\tau)$" for
$\tau <\rho$ and stronger than $q_\rho$, i.e.
$p_{\rho}$  satisfies the following:
for every $\tau <\rho$ there is $s\in
G_\kap *\{\rho\} *Q_{1\kap}\rhookup\rho$ so
that $\langle s, p_{\rho}\rangle \|
\kap\in j(\vta_\tau)$.

Shrink the set $C$ to a set $C^*\in\calU_\lam$
so that for any two $\rho_1\le\rho_2\in C^*$
decisions are the same, i.e. for every $s,\tau
<\rho_1$ as above
$$\eqalignno{&\langle s,p_{\rho_1}\rangle \|\kap\in j
(\vta_{{}_\tau})\cr
\noalign{\hbox{iff}}
&\langle s, p_{\rho_2}\rangle \|\kap\in
j(\vta_{{}_\tau)}\cr
\noalign{\hbox{and}}
&\langle s, p_{\rho_1}\rangle\llvdash\kap\in
j(\vta_{{}_\tau})\cr
\noalign{\hbox{iff}}
&\langle s,p_{\rho_2}\rangle\llvdash\kap\in
j(\vta_{{}_\tau})\ .\cr}$$

For every $\rho\in C^*$  we define
in $V[G_\kap *\{\rho\}*G(Q_{1\kap})]$  a normal
ultrafilter $D(\rho)$  over $\kap$,  where
$G(Q_{1\kap})\subseteq Q_{1\kap}$  generic and
$Q_{1\kap}$  has length $\rho$.  Let us set
$A\in D(\rho)$  iff for some $s\in G_\kap
*\{\rho\} *G(Q_{1\kap})$ $\langle s,p_{\rho}\rangle
\llvdash\kap\in j(\vta_{{}\tau})$  where the
interpretation of $\vta_{{}\tau}$ is $A$ for $\tau
<\rho$, 

Suppose now that $\rho <\rho_1$ are two elements
of $C^*$.  Work in $V[G_\kap *\{\rho_1\}*
G(Q_{1\kap})]$.  Then, clearly,
$G(Q_{1\kap})\rhookup\rho$  will
$V[G_\kap*\{\rho\}]$  generic for $Q_{1\kap}$
(or $Q_{1\kap}\rhookup\rho$ in the sense of the
iteration to $\rho_1$).  So $D(\rho)\in V[G_\kap
*\{ \rho_1\}* G(Q_{1\kap})]$. 

\subheading{Claim} $D(\rho)\subseteq
D(\rho_1)$.

\pr Let $A\in D(\rho)$.  Pick $\tau <\rho ,
\vta_{{}_\tau}$ and $s$ to be as in the definition
of $D(\rho)$.  By the choice of $C^*$, then $\langle
s,p_{\rho_1}\rangle \llvdash \kap\in
j(\vta_{{}\tau})$.  So,
$i_{G(Q_{1\kap})}(\vta_{{}\tau})= i_{G(Q_{1\kap}
\rhookup\rho)} (\vta_{{}\tau})=A$  is in $D(\rho_1)$ as
well as in $D(\rho)$, where $i_G$  is
the function interpreting names.\hfill$\square$ 

Now we are about to complete the construction.
Thus, let $\del$  be a limit of an increasing
sequence $\langle \rho_i\mid i<\kap^+\rangle$  of
elements of $C^*$.  We consider $V[G_\kap
*\{\del\} ]$, i.e. the iteration $Q_{1\kap}$
will be of the length $\del$.  By the claim, 
$$D(\rho_0)\subseteq D(\rho_1)\subseteq
\cdots\subseteq D(\rho_i)\subseteq\cdots
(i<\kap^+)\ .$$
For every $i<\kap^+$, $D(\rho_i)$
is a normal ultrafilter over $\kap$  in
$V[G_\kap *\{\del\} * G(Q_{1\kap})\rhookup
\rho_i]$.  Hence, the forcing $Q_{D(\rho_i)}$
was used at the stage $\rho_i+1$.  Finally set
$D=\bigcup_{i<\kap^+}D(\rho_i)$.

\proclaim Lemma 1.9.  In $V[G_\kap *\{\del\}*
G(Q_{1\kap})]$ $D$  is a normal ultrafilter over
$\kap$  generated by $\kap^+$ sets and
$2^\kap =\del >\kap^+$.

\pr $2^\kap =\del$  since at each stage of the
iteration $Q_{1\kap}$  a new subset of $\kap$
is produced and $\del$  is a limit of inaccessibles
of cofinality $\kap^+$.  

Notice that $D$ is a normal ultrafilter over
$\kap$  since it is an increasing union of
$\kap^+$  normal ultrafilters $D(\rho_i)$
$(D(\rho_i)$ is such in $V[G_\kap *\{\del \}*
Q_{1\kap}\rhookup\rho_i])$ and $Q_{1\kap}$ satisfies
$\kap^+$-c.c.  It is $\kap^+$-generated since for every
$i<\kap^+$  a set $A_{D(\rho_i)}$  generating $D(\rho_i)$
is added at stage $\rho_i+1$.

\hfill$\square$ 

Let $\langle f_{D(\rho_i)}\mid i<\kap^+\rangle$
be the generic functions added by $Q_{D(\rho_i)}$'s.
Use the Prikry forcing with $D$. Let $\langle
\kap_n\mid n<\ome\rangle$  be the Prikry sequence. 
Then by Lemma 1.5 we obtain the following:

\proclaim Theorem 1.10.  The following holds in
the model $V[G_\kap
*\{\del\}*G(Q_{1\kap})*\langle\kap_n\mid
n<\ome\rangle ]$   
\item{(a)} $\kap$  is a strong limit cardinal of
cofinality $\ome$
\smallskip
\item{(b)} $2^\kap =\del >\kap^+$
\smallskip
\item{(c)} the functions $\langle
f_{D(\rho_i)}\mid i<\kap^+\rangle$ are
witnessing $d_{<\aleph_1}(\kap)=\kap^+$.

\subheading{Remarks}  

(1) If one likes to have $2^\kap=\kap^{+7}$
then just collapse $\del$  to $\kap^{+6}$  using
the Levy collapse.  No new subset of $\kap$
will be added.  So $d_{<\aleph_1}$  will still be
$\kap^+$.

(2) $\kap^+$  as the density can be replaced by
$\kap^{++}$, $\kap^{+7}$  etc.  Just pick a
longer sequence of $\rho_i$-s  and argue that no
smaller family is dense.  It requires simple arguments
about names in the Prikry forcing.

(3) $\aleph_1$  can be replaced by any regular
$\tet <\kap$.

\sect{2.~~The Basic Construction}

In this section we will show how to apply
[Git-Mag] in order to produce a model with a
strong limit $\kap$, $\cf\kap =\aleph_0$, 
$2^\kap=\kap^{++}$  and
$d_{<\aleph_1}(\kap)=\kap^+$.

The idea will be that we can reflect the
situation above $\kap$ in the ground model below
$\kap$  once changing its cofinality. 

\subheading{Theorem 2.1}  Suppose that
$V\models GCH$  and there exists an elementary
embedding $j:V\to M$  with a critical point
$\kap$  such that
\item{(a)} $M\supseteq V_{\kap +2}$
\item{(b)} $M=\{ j(f)(\del_1\nek\del_n)\mid n
<\ome$, $\del_1 <\cdots <\del_n<\kap^{++}$  and
$f:[\kap]^n\to V\}$
\item{(c)} ${}^{\!\kap} M\subseteq M$.
\item{}

\smallskip
Then there is a cardinal preserving extension
$V[G]$ of $V$  so that
\item{(1)} for every $\alp <\kap$ or $\alp
>\kap$  $2^\alp =\alp^+$
\smallskip
\item{(2)} $2^\kap =\kap^{++}$
\smallskip
\item{(3)} $\cf\kap =\aleph_0$
\smallskip
\item{(4)} $d_{<\aleph_1}(\kap)=\kap^+$

\subheading{Remark} The assumption used in 2.1
is actually the $\calP^2(\kap)$-hypermeasurability
of $\kap$  or in the Mitchell order $o(\kap)=\kap^{++}+1$.

\pr Let $\calU_0=\{ X\subseteq\kap\mid\kap\in j(X)\}$.
Then $\calU_0$ is a normal ultrafilter over
$\kap$.  Let $i:V\to N\simeq \Ult(V,\calU_0)$ be
the corresponding elementary embedding.  Then
the following diagram is commutative. 
$$\matrix{&&M\cr
&\buildrel j\over\nearrow&\cr
V&&\Big\uparrow\rlap k\cr
&\buildrel i\over\searrow&\cr
&&N\cr}$$
where $k(i(f)(\kap))=j(f)(\kap)$. 

The critical point of $k$ is $\kap^{++}$.

\proclaim Lemma 2.3.  There is a sequence $\l
A_\alp\mid\alp <\kap^+\r$ so that (i)
$j(\kap)=\bigcup_{\alp <\kap^+}A_\alp$  and for
every $\alp <\kap^+$  (ii) $A_\alp\in M$, (iii)
$|A_\alp|\le\kap^{++}$ and (iv) $A_\alp\in
rngk$.

\pr For every $\tau <j(\kap)$  there are
$\del_1\nek\del_n<\kap^{++}$ and $f:[\kap]^n\to\kap$
such that $j(f)(\del_1\nek \del_n)=\tau$.
Consider a function $f':\kap\to [\kap]^{<\kap}$
defined as follows:
$$f'(\nu)=\{ f(\nu_1\nek\nu_n)\mid\nu_1\nek\nu_n<\nu^{++}
\}\ .$$  
Then, in $M$, $|j(f')(\kap)| \le\kap^{++}$  and
$\tau\in j(f')(\kap)$.  Clearly, $k(i(f'))(\kap)=j(f')
(\kap)$.  Hence $j(f')(\kap)\in rngk$.  

So, $j(\kap)=\bigcup \{ j(f')(\kap)\mid f':\kap\to
[\kap]^{<\kap}]$ and for every
$\nu <\kap\ |f'(\nu)|\le\nu^{++}\}$.  Since the
number of such $f'$  is $\kap^+$, we are done.\hfill
$\bigsquare$

\proclaim Lemma 2.4.  There exists a dense set
$F$  of cardinality $\kap^+$ in the topological
space ${}^{j(\kap)}\kap$  with the topology
generated by $<\kap^+$  products such that every
element of $F$ belongs to $rngk$  and in
particular also to $M$.

\pr Let a sequence $\l A_\alp\mid \alp
<\kap^+\r$  be given by Lemma 2.3.  Assume also
that it is increasing.  For every $\alp <\kap^+$
there is $A^*_\alp\in N$,  such that $k(A^*_\alp)=A_\alp$ 
and $N\models |A^*_\alp|=\kap^{++}$.  Working in
$N$  and using GCH, we pick a dense subset
$F^*_\alp$,  with $|F_\alp^*|=\kap^+$  of the
topological space ${}^{A^*_\alp}\kap$  with the
topology generated by countable products.  Then
let $F_\alp =k(F^*_\alp)$  and $F=\bigcup_{\alp
<\kap^+}F_\alp$.  Notice, that $|F_\alp|=\kap^+$
in both $M$  and $V$,  since crit $(k)>\kap^+$.
Clearly, $F$  is as required.\hfill$\bigsquare$

The family $F$  of Lemma 2.4 will be used to
generate a dense set in the space ${}^\kap 2$
with countable product topology once the
cofinality of $\kap$  is changed to $\ome$ and
its power is blown up to $\kap^{++}$.  Thus, if
$\langle \kap_n\mid n<\ome\rangle$  is the
Prikry sequence for the normal measure of the
extender, i.e. for $\calU_0$  and $f=k(f^*)\in
F$,  then let $f^{**}$  be a function such that
$(i)(f^{**})(\kap)=f^*$.  The dense set will
consist of functions $\bigcup_{n<\ome}f^{**}(\kap_n)
\rhookup(\kap_{n+1}\bks\kap_n)$.  

Now, in order to show that this works, we need
to deal with names of clopen sets in ${}^\kap 2$
in the forcing of [Git-Mag].  Finite iterations
described below are needed for a nice
representation of such names.

The model $M$  is the ultrapower of $V$  by the
extender $E=\l E_a\mid a\in [\kap^{++}]^{<\ome}\r$,
where
$$X\in E_a\quad {\rm iff}\quad a\in j(X)\ .$$ 
Now, $j(E)=_{df}E_1\in M=_{df}M_1$  and it is an
extender over $j(\kap^{++})$.  Using $E_1$  we
obtain $j_1:M_1\to M_2\simeq \Ult(M_1, E_1)$
with a critical point $j(\kap)=_{df}\kap_1$.
Let $j_0=j$, $V=M_0$  and $\kap=\kap_0$.  In the
same fashion we can use $j_1(E_1)=_{df}E_2$
over $M_2$ and form $j_2:M_2\to M_3\simeq
\Ult(M_2, E_2)$  with a critical point
$j_1(\kap_1)=_{df}\kap_2$,  an so on.  Thus, for
$n<\ome$, we will have $j_n:M_n\to M_{n+1}\simeq
\Ult (M_n, E_n)$, crit $(j_n)=\kap_n$.  Let
$j_{0n}:V\to M_n$, crit $(j_{0n})=\kap$  be the
composition of $j,\ j_1\nek j_{n-1}$.
Another way to obtain $M_n$'s is using finite
products of $E$ and there ultrapower.  Thus we
consider $E^2=\l E^2_a\mid a\in
[\kap^{++}]^{<\ome}\r$ where for $a\in [\kap]^m$
$(m<\ome)$  and $X\subseteq [\kap]^m\times
[\kap]^m$, $X\in E^2_a$  iff
$\{(\alp_1\nek \alp_m)|\{ (\bet_1\nek \bet_m)|
(\alp_1\nek\alp_m$,  $\bet_1\nek \bet_m)\in X\}
\in E_a\}\in E_a$.  It is not hard to see that
$M_2\simeq \ult (V,E^2)$  and the corresponding
embedding is the same as $j_{02}$.  In the same
fashion for every $n$,  $0<n<\ome$, we can reach
$M_n$  using only one ultrapower.  Thus if
$E^n=\l E^n_a\mid a\in[\kap^{++}]^{<\ome}\r$,
then $M_n\simeq \Ult (V, E^n)$.  Instead of
dealing with finite $a$'s we can replace them
just by ordinals using a reasonable coding. 

The following lemma is routine.

\proclaim Lemma 2.5. For every $\alp
<j_{0n}(\kap)=\kap_n$  there are
$f_\alp:[\kap]^n\to \kap$  and $\del <\kap^{++}$
such that $\alp =j_{0n}(f_\alp)(\del, j_1(\del),
j_2(\del)\nek j_{n-1}(\del))$.

Now fix $n$,  $1<n<\ome$.  We like to describe
one more way of constructing $M_n$.  Thus, we
consider $E^{n-1}$  and $M_1$.  $E^{n-2}$  and
even $E$  is not in $M_1$  but we still can from
outside measure subsets of $\kap$  of $M_1$.  So
we can form $\Ult(M_1, E^{n-1})$.  Since $V_{\kap+2}
\subseteq M_1$  and ${}^\kap\!M_1\subseteq M_1$,
it is routine to check that $\Ult (M_1,
E^{n-1})\simeq M_n$.  Let $\ell$  be the
corresponding embedding. Then $\ell (\kap)=\kap_{n-1}$,
$\ell(\kap_1)=\kap_n$.

\proclaim Lemma 2.6.  For every $\alp <\kap_n$
there are $g_\alp :[\kap]^{n-1}\to\kap_1$  and
$\del <\kap^{++}$  such that $\alp=\ell
(g_\alp)(\del, j_1(\del)\nek j_{n-2}(\del))$.

\pr Let $g_\alp$ be a function representing
$\alp$  in the ultrapower by $E^{n-1}$,  i.e.
for some $\del <\kap^{++}$
$j_{n-1}(g_\alp)(\del, j_1(\del)\nek
j_{n-2}(\del))=\alp$.  Then $g_\alp
:[\kap]^{n-1}\to\kap_1$, since $\alp <\kap_1$
and $j_{n-1}(\kap_1)=\kap_n$.  But then also
$\ell(g_\alp)(\del, j_1(\del)\nek j_{n-2}(\del))=\alp$,
since ${}^\kap\!M={}^\kap\!V$.\hfill$\bigsquare$

Further let us add to such $\ell$  the
subscript $n$.

Let $F$  be the family given by Lemma 2.4.  We
define $F_n=\ell\tagg_n(F)$  for every $n$,  $0<n<\ome$.
Let $\tilF_{nk}=\{ f[\kap_{k-1},\kap_k)\mid
f\in F_n\}$  for every $k$,  $0<k\le n$. For $n$,
$0<n<\ome$  and $t\in \prod^n_{k=1}\tilF_{nk}$
$\ \cup rngt$  is a partial function from $\kap_n$  to
$2$  and it belongs to $M_n$  as a finite union
of its elements.  Set $F^*_n=\{\cup t|$ for some
$m$,  $0<m\le n$  $t\in \prod^m_{k=1}\tilF_{mk}\}$. 

\proclaim Lemma 2.7.  For every $n$, $1<n<\ome$,
$F^*_n$  is dense in the topological space
${}^{\!\kap_n}2$  with countable product topology.

\pr Let $\l \alp_m\mid m<\ome\r$  be an
$\ome$-sequence of ordinals below $\kap_n$, for
some $n$, $1<n<\ome$.  Let $\varphi\in^{\{
\alp_m\mid m<\ome\}}2$.  By the definition of
$F^*_n$  it is enough to prove the lemma in the
situation when all $\alp_m$'s are in some fixed
interval $[\kap_{k-1},\kap_k)$  for $0<k\le n$.
Also by Lemma 2.4, we can assume that $k>1$.
Since nothing happens between $\kap_k$ and
$\kap_n$, we can assume that $k=n$.
For every $m<\ome$, by Lemma 2.6 there are
$g_m:[\kap]^{n-1}\to \kap_1$  and $\del_m<\kap^{++}$
such that $\alp_m=\ell_n(g_m)(\del_m,
j_1(\del_m)\nek j_{n-2}(\del_m))$.  Since
$(\kap^{++})^{\aleph_0}=\kap^{++}$ and $E$ 
is $\ome$-closed.  So we can code the sequence
$\l\del_m\mid m<\ome\r$  into one $\del <\kap^{++}$.
Hence, for every $m<\ome$
$$\alp_m=\ell_n(g_m)(\del ,j_1(\del)\nek
j_{n-2}(\del))\ .$$
Since $\alp_m$'s are all different, there will
be $A\in\calU_\del^{n-1}$  such that for every
$m\not=\ell <\ome$  and $\vec s\in [A]^{n-1}$
$g_m(s)\not= g_\ell(s)$.  Let us also show that
the ranges of $g_m$'s can be made disjoint.  Let
us do this for two $g_0$  and $g_1$.  Using the
completeness of $U_\del$  it is easy then to get
the full result.

\subheading{Claim 1.8.1} There is $B\subseteq A$
in $\calU_\del^{n-1}$  such that $rng$ $g_0\rhookup
B\cap$ $rng$ $g_1\rhookup B=\emptyset$.

\subheading{Remark}  It may not be true iff
either $\alp_0$,  $\alp_1$  are in different
intervals $[\kap_k,\kap_{k+1})$ or if a same 
measure appears in the extender several times.

\pr In order to simplify the notation, let us
assume that $n=3$.  So $\kap_2\le\alp_0<\alp_1<\kap_3$.
Recall that $\ell_3(\kap)=\kap_2$  and
$\ell_3(\kap_1)=\kap_3$.
So, for almost all $(\MOD \calU^2_\del)$
$(\bet,\gam)\in [\kap]^2$  $\kap\le g_0(\bet,\gam)
<g_1(\bet,\gam)<\kap_1$.  Consider $\rho_i=\inf_{C\in
\calU^2_\del}(\sup rng(g_i\rhookup C))$ for
$i<2$.  If $\rho_0\not=\rho_1$, then everything
is trivial. Suppose that $\rho_0=\rho_1=_{df}\rho$.
Then $\cf\rho=\kap$ by $\kap$-completeness of
$\calU_\del$.  Notice also that $g_0$  or $g_1$
cannot be constant $(\MOD\calU^2_\del)$ since
then this constant will be $\rho$.  Consider
sets $X_0=(g_0\tagg [\kap]^2)\cap\rho$  and
$X_1=(g\tagg_1[\kap]^2)\cap\rho$.  We define a
$\kap$-complete ultrafilters $\calW_0$ and
$\calW_1$  over $X_0$  and $X_1$  as follows:   

$\calS\in\calW_k$  iff $g_k^{-1}{}\tagg\calS\in
\calU^2_\del$  where $k<2$.

Then $\calW_0, \calW_1\le_{RK}\calU^2_\del$
(less in the Rudin-Keisler ordering) and $g_0,
g_1$  are the corresponding projection
functions.  Now, $g_0\not= g_1\MOD\calU^2_\del$  
and the extender $E$  has the length $\kap^{++}$.
So, $\calW_0\not=\calW_1$, (see [Mit]). Now we
pick $B_0\in\calW_0\bks\calW_1$  and set $B_1=X_1\bks
B_0$.  The set
$$B=(g_0^{-1}{}\tagg B_0)\cap (g_1^{-1}{}\tagg
B_1)$$
is as desired.\hb
$\bigsquare$ of the claim.

So let $B\in\calU^2_\del$  be so that $g_m\tagg
B\cap g_k\tagg B=\emptyset$ for every $m\not=
k<\ome$.  Denote $g_m\tagg B$  by $B_m$
$(m<\ome)$. Consider now the clopen set in
${}\!^{\kap_1}2$ with $\kap$-products generated
by $\psi\in^{\cup_{m<\ome}B_m}2$  where
$\psi\rhookup B_m$  takes the constant value
$\varphi(\alp_m)$.  Now pick $f\in F_n$ $f\supseteq\psi$.
Then $\ell_n(f)\supseteq\varphi$,  since for every
$m<\ome$  $\{ (\bet,\gam)\in [\kap]^2|f(g_m(\bet,\gam))
=\varphi(\alp_m)\}\supseteq g_m^{-1}$"$(B_m)
\in\calU^2_\del$.\hfill$\bigsquare$

Suppose now that the extender $E$  has the
length $\kap^{+7}$  instead of $\kap^{++}$.  We
would like to apply previous arguments in order
to produce a dense set of cardinality $\kap^{+6}$.
The only obstacle is that Lemma 2.7 breaks down
if we use the family $F^*_n$  defined there.
The problem is that once the length of the
extender is $\ge\kap^{+++}$ same measures are
starting to appear in it at different places.
It was crucial for Claim 2.8.1 that this does doesn't
happen.  The solution is going to be to take a
larger family and use the fact that for any two
measures in the extenders there is a measure
with index $<\kap^{+6}$ which provides a
difference between them inside the extender. 

First let us define the new $F_n$.  Let
$F'_n=\{t\mid t:[\kap]^n\to F\}$.  Clearly,
$|F_n'|=|F|^\kap=(\kap^{+6})^\kap =\kap^{+6}$.

Also, every $t$  in $F'_n$  is in $M$ since
${}\!^\kap M\subseteq M$.  Now for every $\del
<\kap^{+6}$  and $t\in F'_n$  we consider
$\ell_n(t)(\del, j_1(\del)\nek j_{n-2}(\del))$.
It is an element of $M_n$.  Set
$$F_n=\{ \ell_n(t)(\del, j_1(\del)\nek j_{n-2}(\del))\mid
\del <\kap^{+6}\quad {\rm and}\quad t\in F'_n\}\ .$$  

Now, we define $F^*_n$ as in case $\kap^{++}$
using this new $F_n$.  We need to show that the
analog of 2.7 is true with our new $F_n$.  The
arguments of 2.7 and 2.8.1 are running smooth
until the point where it is claimed that
$W_0\not= W_1$.  

Suppose now that $W_0=W_1$. Let us assume in
order to simplify the presentation that 
$$\kap <\kap_1<\alp_0,\alp_1<\kap_2=\kap_n\ .$$
Thus $g_0,g_1$  are now one-place functions.

The ultrafilters $W_0, W_1$  are then isomorphic
to some measures $\calU_{\tau_0}, \calU_{\tau_1}$
of extender, where $\tau_0, \tau_1<\kap^{+7}$
and for $\tau <\kap^{+7}\ \calU_\tau
=\{\calS\subseteq \kap\mid \tau\in j(\calS)\}$.
Just take the bijections $\rho_0,\rho_1$ between
$\kap$  and $X_0,X_1$.  The general case is
slightly more complicated since we need to deal
with $E^2$,  $E^3$ etc. and instead of $\calU_\tau$
it will be $\calU_\tau^2,\calU_\tau^3$, etc.

Let $\tau_0<\tau_1<\kap^{+7}$.  The opposite
case is identical.  There exists $\tau
<\kap^{+6}$  such that $E_{\{\tau,\tau_0\}}\not=
E_{\{\tau,\tau_1\}}$,  where for $a\in
[\kap]^{<\ome}$ $E_a=\{ \calS\subseteq\kap^{|a|}\mid
a\in j(\calS)\}$.  For this use for example
$\tau$  coding the walk from $\tau_1$  to
$\tau_0$  since the coherent box sequence of
$\kap^{+7}$.

Next, we replace $\del$ by some $\del^*<\kap^{+7}$
coding $\{\tau,\tau_0,\tau_1,\del\}$.  Or in
other words, we find $U_{\del^*}$  in the
extender $E$  which is Rudin-Keisler above $\calU_\del,
E_{\{\tau,\tau_0\}}, E_{\{\tau, \tau_1\}}$.  Let
$\pi_\del$  be the corresponding projection of
$\calU_{\del^*}$  onto $U_\del$.  Define
$g^*_i:\kap\to\kap_1$ $(i<2)$  as follows 
$$g_i^*(\bet)=g_i(\pi_\del(\bet))\ .$$
Then, $\alp_i=\ell_2(g_i^*)(\del^*)$. Hence
$g^*_i$  projects $\calU_{\del^*}$  onto $W_i$.
Consider an ultrafilter $E_i$  over $\kap\times
X_i$  defined as follows:

$\calS\in E_i$  iff for some
$\calS'\in E_{\{\tau,\tau_i\}}$ $\calS =\{
(\bet, \rho_i(\gam))\mid
(\bet,\gam)\in\calS'\}$.

I.e. we are using the
bijection $\rho_i$  to transfer $E_{\{
\tau_i\}}$  back to $W_i$.
Pick projections $\pi_i$  and $\pi$ of $\calU_{\del^*}$
to $E_i$  and $\calU_\tau$  such that
$\pi_i(\xi)=(\pi(\xi), g^*_i(\xi))$ for almost
all $\xi$ $\MOD\calU_{\del^*}$.

Now we find disjoint $B'_0\in E_0$  and
$B'_1\in E_1$.  There is $B\in\calU_{\del^*}$
such that $\pi_0\tagg(B)\subseteq B'_0$  and
$\pi_1\tagg(B)\subseteq B'_1$.  Let
$C=\pi\tagg(B)$,  $B_0=\pi\tagg_0(B)$  and
$B_1=\pi\tagg_1(B)$.  Then $C\in\calU_\tau$,
$B_0\in E_0$  and $B_1\in E_1$. The following is
important:\hb
$(*)$  for every $\bet\in C$  and $\gam <\kap_1$
it is impossible to have both $(\bet,\gam)\in
B_0$  and $(\bet,\gam)\in B_1$.  For $\bet\in C$
we consider the set $C_\bet =\{ \gam\in X_0\cup
X_1\mid (\bet,\gam)\in B_0\cup B_1\}$. For every
$\bet\in C$  let $\psi_\bet:C_\bet\to 2$  be
defined as follows: 
$$\psi_\bet(\gam)=\cases{\varphi(\alp_0)&if
$(\bet,\gam)\in B_0$\cr
\varphi(\alp_1)&if $(\bet,\gam)\in B_1$\cr}$$
Notice that by $(*)$  such defined $\psi_\bet$
is a function.  Since $|C_\bet|\le\kap$,
$C_\bet\in M$  and $C_\bet\subseteq\kap_1$,
there is $f_\bet\in F$
$f_\bet\supseteq\psi_\bet$.  Let $t:\kap\to F$
be defined by $t(\bet)=f_\bet$  for $\bet\in C$
and arbitrarily (but in $F$) otherwise.  Then,
$\ell_2(t)(\tau)\in F_2$  and let us show that
$\ell_2(t) (\tau)\supseteq\varphi\rhookup\{\alp_0,
\alp_1\}$.  It is enough to show that the set   
$$\{\xi <\kap\mid \pi(\xi)\in C\ g^*_0(\xi),
g^*_1(\xi)\in C_{\pi(\xi)}\quad {\rm and}\quad  
f_{\pi(\xi)}(g^*_i(\xi))=\varphi(\alp_i)\quad {\rm for}
\quad i<2\}$$
is in $\calU_{\del^*}$.  We claim that it contains $B$. 
Thus let $\xi\in B$.  Then, $\pi (\xi)\in C$,  $(\pi
(\xi), g^*_0(\xi))\in B_0$ and $(\pi(\xi), g^*_1(\xi))\in
B_1$.  Hence, $g^*_0(\xi), g^*_1(\xi)\in C_{\pi (\xi)}$
and $f_{\pi(\xi)}$ was chosen so that $f_{\pi(\xi)},
(g_i(\xi))=\varphi(\alp_i)$  where $i<2$. 

This show the density for $\alp_0,\alp_1$.  In
order to deal with $\l\alp_m\mid m<\ome\r$
instead of only two $\alp_0,\alp_1$,  just
produce disjoint $\l B_m\mid m<\ome\r$  using
$\ome_1$-completeness of the ultrafilters involved. 

Now we are ready to complete the
proof of Theorem 2.1. For every $n$,  $0<n<\ome$
let $F^*_n$  be a set given by Lemma 2.7. Then
for every $f\in F^*_n$  $(1<n<\ome)$ there will
be $\of:[\kap]^n\to {}^{\!\kap\ge}\kap$
representing $\of$  in the ultrapower by
$\calU^n_0$, i.e.
$$j_{n-1}(\of)(\kap_0,\kap_1\nek \kap_{n-1})=f\ .$$
Set $\oF_n=\{ \of\mid\of:[\kap]^n\to{}^{\kap\ge}\kap$
and $j_{n-1}(\of)(\kap_0,\kap_1\nek
\kap_{n-1})\in F^*_n\}$, where $0<n<\ome$.  Let
$\oF_0=\{\of \mid\dom\of =\{ 0\}$  and $\of
(0):\kap\to\kap\}$.  Define
$\oF=\prod_{n<\ome}\oF_n$. Clearly,
$|\oF |=\prodl_{n<\ome}|\oF_n|=(\kap^+)^{\aleph_0}
=\kap^+$.

Suppose now that we are forced with the forcing
of [Git-Mag], then basically, a Prikry sequence
was added for every measure of the extender $E$
and no new bounded subset of $\kap$  was added.
So, GCH holds below $\kap$,  $\cf\kap=\aleph_0$
and $2^\kap=\kap^{++}$.  Let $\l\nu_n \mid
n<\ome\r$  be the Prikry sequence for $\calU_0$,
i.e. for the normal measure.  We are going to
use it in order to define a dense set $D$ in the
topological space ${}^{\!\kap}2$ with topology
generated by countable products.  The idea is to
transfer $F^*_n$'s to the space ${}^{\!\kap}2$.
We are going to take functions representing
elements of $F^*_n$'s i.e. the members of
$\oF_n$  and apply them to the $\l\nu_1\nek\nu_n\r$.
Then, in order to show density we will notice
that a name of a basic clopen set can be
transferred back to $\kap_n$'s using the same
process but in the opposite direction. Over $\kap_n$
we find an element of $F^*_n$  inside such
clopen set and pull it back to ${}^{\!\kap}2$.

Now let us do this formally.  For every $t\in\oF$
we define a partial function $t^*$  from $\kap$
to $2$  as follows.  Let $\alp <\kap$.  We find
$n_\alp <\ome$  such that $\nu_{n_\alp}\le \alp
<\nu_{n_\alp+1}$,  where $\nu_0$  denotes $0$.
If $n_\alp =0$ and $t(0)(\alp)<\nu_1$, then set
$\alp\in\dom t^*$  and $t^*(\alp)=t(0)(\alp)$.
Suppose now that $n_\alp >0$.  If $\alp\in\dom
t(n_\alp)(\nu_1\nek \nu_{n_\alp})$  and
$\nu_{n_\alp}\le t(n_\alp)(\nu_1\nek \nu_{n_\alp})
<\nu_{n_\alp +1}$ then set $\alp\in\dom t^*$
and $t^*(\alp)=t(n_\alp)(\nu_1\nek
\nu_{n_\alp})(\alp)$.  Otherwise $t^*(\alp)$  is
undefined or if one likes to have it total just set then
$t^*(\alp)=0$.  Set $D=\{ t^*\mid
t\in\oF\}$.  Obviously, $|D|\le |\oF|=\kap^+$.  

\proclaim Lemma 2.9.  $D$  is dense in the
topological space ${}^{\!\kap}2$ with the
topology generated by countable products.

\pr Suppose $\varphi\in^{\{\tau_m\mid m<\ome\}}2$.
We need to find some $f\in D$ $f\supseteq\varphi$.
Let us work in $V$ with names instead of working
in the generic extension.  So, let
$\vtop{\offinterlineskip\hbox{$\tau_m$}\hbox{$\sim$}}$
be a name of an ordinal $\tau_m(m<\ome)$
and $\vtop{\offinterlineskip\hbox{$\varphi$}
\hbox{$\sim$}}$  a name for $\varphi$.

Our basic tool will be Lemma 2.11 of [Git-Mag]
or actually the condition $p^*=p\cup\{\l
\bet,\emptyset, S^*\r\}$  produced in this lemma
if instead of $\buildrul\sim\under g$  we deal
with $\l\vtop{\offinterlineskip\hbox{$\tau_m$}\hbox
{$\sim$}}\mid m<\ome\r$
and $\buildrul\sim\under\varphi$  there.  In
order to make the presentation as self-contained
as possible, let us state here the main
properties of $p^*$.  Thus $S^*$ is a subtree of
$[\kap]^{<\ome}$  such that for every $s\in S^*$
$Suc_{S^*}(s)\in\calU_\bet$.  For every $m<\ome$
there is a level $n_m<\ome$  in $S^*$ such that
for every $s_1,s_2\in S^*$  from this level,
i.e. $|s_1|=|s_2|=n_m$  there are
$\gam_1,\gam_2$  and $i<2$  such that the
following holds for $k=1,2$
$$(s_k(n_m))^0\le\gam_k<(\min Suc_{S^*}(s_k))^0\leqno(a)$$
$$(p\cup\{\l \bet,\emptyset, S^*\r\})_{s_k}\llvdash
(\vtop{\offinterlineskip\hbox{$\alp_m$}\hbox
{$\sim$}}=\gam_k\quad {\rm and}\quad\buildrul\sim\under
\varphi(\vtop{\offinterlineskip\hbox{$\alp_m$}\hbox
{$\sim$}}=i)\ ,\leqno(b)$$
where ${}^0$-denotes the projection function to
the normal measure $\calU_0$  and $(p\cup\{\l
\bet,\emptyset ,S^*\r \})_{s_k}$  is the
condition obtained from $p\cup \{\l \bet,
\emptyset, S^*\r \}$  by adding $s_k$  to be the
initial segment of the Prikry sequence for $\bet$
(or $\calU_\bet$) and then shrinking $S^*$  to
the tree above $s_k$  and projecting $s_k$  to
the appropriate coordinates in $p$.

Now consider the following set
$$A=\{ n<\ome\mid\exists m<\ome\quad n=n_m\}$$
Let $n\in A$.  Denote $\{m <\ome \mid n_m=n\}$
by $A_n$.  We define a function $g_n$  on
$Lev_n(S^*)$. Let $s\in Lev_n(S^*)$.  By (a),
(b), for every $m\in A_n$  there are
$\gam_{m,s}$  and $i_m<2$  such that 
$$(s(n))^0\le\gam_m<(\min Suc_{S^*}(s))^0\leqno(1)$$
$$(p\cup\{ \l \bet, \emptyset, S^*\r
\})_s\llvdash(\vtop{\offinterlineskip\hbox{$\alp_m$}
\hbox{$\sim$}}=\gam_{m,s}\quad {\rm
and}\quad\buildrul\sim\under\varphi
(\vtop{\offinterlineskip\hbox{$\alp_m$}\hbox{$\sim$}})=
i_m)\leqno(2)$$

Set $g_n(s)=\{\l\gam_{m,s},i_m\r\mid m\in
A_n\}$.  Hence, $g_n(s)\in^{\{\gam_{m,s}\mid m\in
A_n\}} 2$.  Then, $g_n$  represents a basic
clopen set in ${}^{\kap_n\!}2$ in $M_n$.
Namely, $j_{n-1}(g_n)(\bet, j_1(\bet)\nek
j_{n-1}(\bet))$.  Using the density of $F^*_n$,
we find $f_n\in F^*_n$  $f_n\supseteq
j_{n-1}(g)(\bet, j_1(\bet)\nek j_{n-1}(\bet))$.
Pick $\of_n\in \oF_n$ such that    
$j_{n-1}(\of_n)(\kap_0,\kap_1\nek\kap_{n-1})=f_n$.
Then for almost all ($\MOD\calU^n_\bet)s\in
Lev_n(S^*)\ \of_n((s)^0)\supseteq g_n(s)$.

Now let us do it for every $n\in A$  we will get
a sequence $\l\of_n\mid n\in A\rangle$.  Let
$t\in\oF$ be such that for every $n\in A$ $t(n)=\of_n$.
Then the corresponding $t^*$ or here its name
$\vtop{\offinterlineskip\hbox{$t^*$}\hbox{$\sim$}}$
will be as desired,
i.e. $p\cup \{\l\bet,\emptyset, S^*\r\}$
$\llvdash \vtop{\offinterlineskip\hbox{$t^*$}
\hbox{$\sim$}}\supseteq\buildrul\sim\under\varphi$.  This
completes the proof of the lemma and hence of
the theorem.\hfill$\bigsquare$

\sect{3.~~Some Generalizations}

Under the same lines we obtain the following
theorem:

\proclaim Theorem 3.1.  Suppose that $o(\kap)=\lam^++1$
(i.e. extender of the length $\lam^+$) and $\cf
\lam>\kap$.  Then the following holds in a
generic extension $V[G]$:
\item{(1)} for every $\alp <\kap$  or
$\alp\ge\lam\ 2^\alp =\alp^+$.
\smallskip
\item{(2)} $2^\kap =\lam^+$
\smallskip
\item{(3)} $\cf\kap =\aleph_0$
\smallskip
\item{(4)} $d_{<\aleph_1}(\kap)=\lam$.

\pr Apply the construction of Section 2 with
extender $E$  of the length $\lam^+$  instead of
$\kap^{++}$.  An additional property that we
need to show in the present situation is that
$d_{<\aleph_1}(\kap)$  cannot be below $\lam$.
But this follows by [Sh430, 5.3, 5.4] and the
$pcf$ structure of the models of [Git-Mag] or
just directly using the correspondence
established in Lemma 2.9 between basic clopen
sets of ${}^\kap{\! 2}$ of $V[G]$ and ${}^{\kap_n}\!2$
of $M_n$.  Since already ${}^{\kap_1}\!2$ cannot have a
dense set of cardinality less than $\lam$ because
${}^{\lam^+}\!2$  embeds it and GCH holds.
\hfill$\bigsquare$

The following two results are straightforward
applications of the techniques for pushing
everything down to $\aleph_\ome$ [Git-Mag,
Section 2] or changing cofinality to $\aleph_1$
Segal [Seg], [Git-Mag2] and pushing down to
$\aleph_{\ome_1}$. 

\proclaim Theorem 3.2.  Suppose
$o(\kap)=\kap^{++}+1$.  Then the following holds
in a generic extension: 
\item{(1)} for every $\alp <\ome$  or $\alp
>\ome$  $2^{\aleph_\alp}=\aleph_{\alp +1}$
\smallskip
\item{(2)} $2^{\aleph_\ome}=\aleph_{\ome +2}$
\smallskip
\item{(3)} $d_{<\aleph_1}(\aleph_\ome)=\aleph_{\ome+1}$.

\proclaim Theorem 3.3.  Suppose $o(\kap)=\kap^{++}+
\ome_1$.  Then the following holds in a generic
extension:  
\item{(1)} for every $\alp$  $2^{\aleph_{\alp +1}}=
\aleph_{\alp +2}$  
\smallskip
\item{(2)} GCH above $\aleph_{\ome_1+1}$
\smallskip
\item{(3)} $2^{\aleph_{\ome_1}}=\aleph_{\ome_1+2}$.
\smallskip
\item{(4)} $d_{<\aleph_2}(\aleph_{\ome_1})
=\aleph_{\ome_1+1}$.

For Theorem 3.3 we need also to replace
$\aleph_0$-box products by $\aleph_1$-base
products.  Notice that all the considerations of
Section 2 are going smoothly if we replace
$\aleph_0$-box product by $\tet$-box product for
any $\tet <\kap$. Also instead of the space
${}^{\kap}\!2$ we can work with ${}^\kap\!\chi$
for any fixed $\chi <\kap$.  So the following
holds: 

\proclaim Theorem 3.4. Suppose that
$o(\kap)=\kap^{++}+1$, $\tet$,  $\chi <\kap$.  Then the
following holds in a generic cardinal preserving
extension:
\item{(1)} for every $\alp <\kap$  or $\alp
>\kap$  $2^\alp =\alp^+$
\smallskip
\item{(2)} $\cf \kap =\aleph_0$
\smallskip
\item{(3)} $2^\kap =\kap^{++}$
\smallskip
\item{(4)} {\sl the density of the topological space
${}^\kap\!\chi$ with the topology generated by
$\tet$-products is $\kap^+$.}

The analogs of 3.2 and 3.3 hold as well.

\sect{4.~~Reaching the Maximal Density and Wider
Gaps}

In previous sections, we constructed models with
density less than the maximal possible value
$2^\kap$.  Let us show now how to construct a
model with the density $2^\kap$  assuming
singularity of $\kap$  and $2^\kap >\kap^+$. 

\proclaim Theorem 4.1. Suppose $o(\kap)=\kap^{+3}+1$,
then there is a generic extension $V[G]$ satisfying the
following 
\item{(1)} for every $\alp <\kap$  or $\alp >\kap$ 
$2^\alp =\alp^+$
\smallskip
\item{(2)} $\cf \kap =\aleph_0$
\smallskip
\item{(3)} $2^\kap=\kap^{++}$
\smallskip
\item{(4)} $d_{<\aleph_1}(\kap)=2^\kap$

\pr Let $V_1$  be a model of Theorem 3.1 with
$\lam=\kap^{++}$.  Let $D$  be a set witnessing
$d_{<\aleph_1}(\kap)=\kap^{++}$.  Collapse $\kap^{+++}$
to $\kap^{++}$  using the Levy collapse.  Let
$V_2$  be such generic extension.  Then, in
$V_2$,  $2^\kap =\kap^{++}$  and
$|D|=\kap^{++}$.  However, $D$  is still
witnessing $d_{<\aleph_1}(\kap)=\kap^{++}$.
Thus, no new subset of $\kap$ are added.  Hence   
$({}\!^\kap\!2)^{V_1}=({}\!^\kap\!2)^{V_2}$.
But also no new subsets of cardinality $\kap^+$
are added to sets of $V_1$.  So there is no
dense set in ${}^\kap\!2$  of cardinality
$\le\kap^+$.  $D$ is dense since there is no new
basic clopen sets.\hfill$\bigsquare$

As in Section 3 it is possible to push this result
down to $\aleph_\ome$  and $\aleph_{\ome_1}$.

Suppose now that one likes to have $2^\kap$ big
but still keep the density $\kap^+$.  A slight
modification of the construction of Section 2
will give the following:

\proclaim Theorem 4.2.  Suppose that $\lam
>\kap$  is a regular cardinal $o(\kap)=\lam +1$.
Then there is a generic cardinal preserving
extension satisfying the following:
\item{(1)} $\kap$  is a strong limit 
\smallskip
\item{(2)} $\cf\kap =\aleph_0$
\smallskip
\item{(3)} $2^\kap =\lam$
\smallskip
\item{(4)} $d_{<\aleph_0}(\kap)=\kap^+$

\pr Let $V\models GCH$.  $E$ an extender of the
length $\lam$,  $j:V\to M\simeq Ult(V,E)$.
Using Backward Easton forcing we blow up $2^{\kap^+}$
to $\lam$.  By standard arguments $E$ extends to
an extender $E^*$  in such generic extension
$V[G]$ as well as $j\subseteq j^*:V[G]\to
M[G^*]$.  Now we proceed with $V[G], M[G^*]$
and $j^*$  as in Section 1.  $\lam$  generic
functions from $\kap^+$  to $\kap^+$  are used 
also to show that the analog of Claim 2.8 is
valid.\hfill$\bigsquare$
\vfill\eject
\references {60}

\ref{[CEG]} F. Carter, P. Erd\"os and F. Galvin, On
density of $\lam$-box products, General Topology
and its Applications 9 (1978), 307-312.

\smallskip
\ref{[CR]} W. Comfort and C. Robertson,
Cardinality constraints for pseudo-compact and
totally dense subgroups of compact topological
groups, Pacific Journal of Mathematics 119 (1985),
265-285.

\smallskip
\ref{[Git-Mag]} M. Gitik and M. Magidor, The
singular cardinal hypothesis revisited, in Set
Theory of the Continuum, H. Judah, W. Just and
H. Woodin eds., (1992), 243-279.

\smallskip
\ref{[Git-Mag2]} M. Gitik and M. Magidor, The
Extender based forcings, Journal of Sym. Logic.,
v.59 (1994), 445-460. 

\smallskip
\ref{[Ka]} A. Kanamari, The Higher Infinite,
Springer, 1995.

\smallskip
\ref{[Sh80]} S. Shelah, A weak generalization of
MA to higher cardinals, Is. J. of Math. 30,
(1978), 297-306. 

\smallskip
\ref{[Mit]} W. Mitchell, Hypermeasurable Cardinals.

\smallskip
\ref{[Seg]} M. Segal, M.Sc. thesis, Jerusalem
1993.

\smallskip
\ref{[Sh430]} S. Shelah, Further Cardinal
Arithmetic, Israel Journal of Math. 

\end

%% file: def1200.tex
%


\def\today{\ifcase\month\or January\or February\or
March\or April\or May\or June\or July\or August\or
September\or October\or November\or December\fi
\space\number\day, \number\year}




\def\dspace{\lineskip=2pt\baselineskip=18pt
\lineskiplimit=0pt}
\def\sspace{\lineskip=2pt\baselineskip=12pt
\lineskiplimit=0pt}

\font \bbrm=cmbx10 at 12pt

\font \ninerm= cmr10 at 9pt

\def\smalltype{\ninerm}

\def\bigtype{\bbrm}

\hsize=13.5cm
\magnification=1200
\def\ce{\centerline}

\def\hb{\hfill\break}

\def\title #1{\null\bigskip\ce{\bigtype #1}
\bigskip}

\def\alp{\alpha}		
\def\bet{\beta}		
\def\gam{\gamma}		
\def\del{\delta}		
\def\eps{\varepsilon}		

\def\tet{\theta}		

\def\kap{\kappa}
\def\lam{\lambda}		
\def\sig{\sigma}		

\def\ome{\omega}		


\def\calP{{\cal P}}

\def\calS{{\cal S}}

\def\calU{{\cal U}}

\def\calW{{\cal W}}



    
\font\tenboldgreek=cmmib10
 \font\sevenboldgreek=cmmib10 at 7pt
\font\fiveboldgreek=cmmib10 at 7pt
\newfam\bgfam
\textfont\bgfam=\tenboldgreek
\scriptfont\bgfam=\sevenboldgreek
\scriptscriptfont\bgfam=\fiveboldgreek

\mathchardef\ggarrow="7010
\def\twoheadrightarrow{{\fam=\myfam\ggarrow}}

\font\tengerman=eufm10 \font\sevengerman=eufm7
\font\fivegerman=eufm5
\font\tendouble=msym10 \font\sevendouble=msym7
\font\fivedouble=msym5

\textfont4=\tengerman \scriptfont4=\sevengerman
\scriptscriptfont4=\fivegerman
\newfam\dbfam
\textfont\dbfam=\tendouble \scriptfont\dbfam=
\sevendouble
\scriptscriptfont\dbfam=\fivedouble

\mathchardef\ng="702D
\mathchardef\dbA="7041
\mathchardef\sm="7072
\mathchardef\nvdash="7030
\mathchardef\nldash="7031
\mathchardef\lne="7008
\mathchardef\sneq="7024
\mathchardef\spneq="7025
\mathchardef\sne="7028
\mathchardef\spne="7029
\mathchardef\ltms="706E
\mathchardef\tmsl="706F

\mathchardef\dbA="7041
\def\ltimes{{\fam=\dbfam\ltms}}

\def\subsetneqq{\,{\fam=\dbfam\sneq}\,}
\def\supsetneqq{\,{\fam=\dbfam\spneq}\,}


\mathchardef\dbA="7041 
\mathchardef\dbB="7042 
\mathchardef\dbC="7043 
\mathchardef\dbD="7044 
\mathchardef\dbE="7045 
\mathchardef\dbF="7046 
\mathchardef\dbG="7047 
\mathchardef\dbH="7048 
\mathchardef\dbI="7049 
\mathchardef\dbJ="704A 
\mathchardef\dbK="704B 
\mathchardef\dbL="704C 
\mathchardef\dbM="704D 
\mathchardef\dbN="704E 
\mathchardef\dbO="704F 
\mathchardef\dbP="7050 
\mathchardef\dbQ="7051 
\mathchardef\dbR="7052 
\mathchardef\dbS="7053 
\mathchardef\dbT="7054 
\mathchardef\dbU="7055 
\mathchardef\dbV="7056 
\mathchardef\dbW="7057 
\mathchardef\dbX="7058 
\mathchardef\dbY="7059 
\mathchardef\dbZ="705A 

\def\nek{,\ldots,}
\def\sdp{\times \hskip -0.3em {\raise 0.3ex
\hbox{$\scriptscriptstyle |$}}} 


\def\cf{{\rm \,cf\,}}

\def\dom{\mathop{\rm dom}\nolimits}

\def\min{\mathop{\rm min}}
\def\MOD{\mathop{\rm mod}}



\def\of{{\overline f}}
\def\oF{{\overline F}}







\def\tilF{{\widetilde F}}


\def\ddownarrow{\big\downarrow \hskip-0.70em\raise
2pt\hbox {$\big\downarrow$}}
\def\longright #1#2 {\smash{\mathop{\hbox to
#1pt {\rightarrowfill}}\limits_{#2}}}
\def\sqr#1#2{{\vcenter{\hrule height.#2pt\hbox{\vrule
width.#2pt height#1pt \kern#1pt \vrule width.#2pt}
\hrule height.#2pt}}}
\def\square{\mathchoice{\sqr34}{\sqr34}{\sqr{2.1}3}
{\sqr{1.5}3}}

\def\buildrul#1\under#2{\mathrel{\mathop{\null#2}
\limits_{#1}}}

\def\boxit#1{\vbox{\hrule\hbox{\vrule\kern3pt
\vbox{\kern3pt#1 \kern3pt}\kern3pt\vrule}\hrule}}

\def\prodl{\prod\limits}

\def\subsetneq{\baselineskip=5pt\vcenter
{\hbox{$\subset$}\hbox{$\ne$}}}

\def\subheading#1{\medskip\goodbreak\noindent{\bf
#1.}\quad}

\def\sect#1{\goodbreak\bigskip\centerline{\bf#1}
\medskip}
\def\pr{\smallskip\noindent{\bf Proof:\quad}}
\def\onumber #1{\ooalign{\hfil\raise.07ex\hbox{
\hfill$\scriptstyle \,#1$\hfil}
\cr\cr{$\bigcirc$}}}
\def\onumber c{\ooalign{\hfil\raise.07ex\hbox
{\hfill$\scriptstyle \,c$\hfil}
\cr\cr{$\bigcirc$}}}
\def\alpcirc {\ooalign{\hfil\raise.07ex
\hbox{\hfill$\scriptstyle\alp\;$\hfill}\cr\cr
{$\bigcirc$}}}
\def\astcirc {\ooalign{\hfil\raise.07ex
\hbox{\hfill$\textstyle\ast\;$\hfill}\cr\cr
{$\bigcirc$}}}

\def\longmapright #1#2 {\smash{\mathop{\hbox to
#1pt {\rightarrowfill}}\limits^{#2}}}
\def\longmapleft #1 #2 {\smash{\mathop{\hbox to
#1 pt {\leftarrowfill}}\limits^{#2}}}

\def\references#1{\goodbreak\bigskip\par\centerline
{\bf References}\medskip\parindent=#1pt}
\def\ref#1{\par\smallskip\hang\indent\llap{\hbox
to \parindent{#1\hfil\enspace}}\ignorespaces}

\def\back{{\raise 2.5pt\hbox{$\,\scriptscriptstyle
\backslash\,$}}}
\def\bks{{\backslash}}
\def\part{\partial}
\def\lwr #1{\lower 5pt\hbox{$#1$}\hskip -3pt}
\def\rse #1{\hskip -3pt\raise 5pt\hbox{$#1$}}
\def\lwrs #1{\lower 4pt\hbox{$\scriptstyle #1$}
\hskip -2pt}
\def\rses #1{\hskip -2pt\raise 3pt\hbox
{$\scriptstyle #1$}}

\def\<#1{\left\langle{#1}\right\rangle}

\def\subinbn{{\subset\hskip-8pt\raise 0.95pt
\hbox{$\scriptscriptstyle\subset$}}}

\def\llvdash{\mathop{\|\hskip-2pt
\raise 3pt\hbox{\vrule height 0.25pt width 1.5cm}}}

\def\lvdash{\mathop{|\hskip-2pt \raise 3pt\hbox
{\vrule height 0.25pt width 1.5cm}}}

\def\fakebold#1{\leavevmode\setbox0=\hbox{#1}%
  \kern-.025em\copy0 \kern-\wd0
  \kern .025em\copy0 \kern-\wd0
  \kern-.025em\raise.0333em\box0 }

\font\msxmten=msxm10
\font\msxmseven=msxm7
\font\msxmfive=msxm5
\newfam\myfam
\textfont\myfam=\msxmten
\scriptfont\myfam=\msxmseven
\scriptscriptfont\myfam=\msxmfive
\mathchardef\rhookupone="7016
\mathchardef\ldh="700D
\mathchardef\leg="7053
\mathchardef\ANG="705E
\mathchardef\lcu="7070
\mathchardef\rcu="7071
\mathchardef\leseq="7035
\mathchardef\qeeg="703D
\mathchardef\qeel="7036
\mathchardef\blackbox="7004
\mathchardef\bbx="7003
\mathchardef\simsucc="7025

\def\rhookup{{\fam=\myfam \rhookupone}}
\def\succsim{\mathrel{\fam=\myfam\simsucc}}

\def\bigsquare{{\fam=\myfam\bbx}}

\font\tencaps=cmcsc10
\def\smallcaps{\tencaps}

\def\author#1{\bigskip\ce{\smallcaps #1}\medskip}

\def\tagg{^{\prime\prime}}

\def\upddots{\mathinner{\mkern
1mu\raise 1pt \hbox{.}\mkern 2mu \mkern
2mu \raise 4pt\hbox{.}\mkern 1mu \raise 7pt\vbox
{\kern 7 pt\hbox{.}}} }

\def\varchi{\ooalign{{\raise
1.385pt\hbox{$\chi$}}\crcr\hbox{--}\crcr}}

\def\trianarrow{{\raise 2pt\hbox to 0.50cm
{\hrulefill}\triangleright}}
\def\Chi{{\raise 3pt\hbox{$\chi$}}}

\font\b=cmr10 scaled \magstep4

\def\bigzerou{\smash{\lower1.7ex\hbox{\b 0}}}
\def\bigast{\smash{\lower1.7ex\hbox{\b *}}}

\def\leaderfill{leaders\hbox to 5em{\hss\hss}\hfill}
\newcount\notenumber
       
       \def\note#1{\advance\notenumber by
1\footnote{$^{\the\notenumber}$} {\sspace\smalltype #1}}

%% file: mssymb.tex
%

%

\catcode`\@=11
\def\relaxnext@{\let\next\relax}
\def\noaccents@{\def\accentfam@{0}}

\font\tenmsx=msxm10
\font\sevenmsx=msxm7
\font\fivemsx=msxm5
\font\tenmsy=msym10
\font\sevenmsy=msym7
\font\fivemsy=msym5
\newfam\msxfam
\newfam\msyfam
\textfont\msxfam=\tenmsx  \scriptfont\msxfam=\sevenmsx
  \scriptscriptfont\msxfam=\fivemsx
\textfont\msyfam=\tenmsy  \scriptfont\msyfam=\sevenmsy
  \scriptscriptfont\msyfam=\fivemsy

\def\hexnumber@#1{\ifcase#1 0\or1\or2\or3\or4\or5\or6\or7\or8\or9\or
	A\or B\or C\or D\or E\or F\fi }

\font\teneuf=eufm10
\font\seveneuf=eufm7
\font\fiveeuf=eufm5
\newfam\euffam
\textfont\euffam=\teneuf
\scriptfont\euffam=\seveneuf
\scriptscriptfont\euffam=\fiveeuf
\def\frak{\relaxnext@\ifmmode\let\next\frak@\else
 \def\next{\Err@{Use \string\frak\space only in math mode}}\fi\next}
\def\goth{\relaxnext@\ifmmode\let\next\frak@\else
 \def\next{\Err@{Use \string\goth\space only in math mode}}\fi\next}
\def\frak@#1{{\frak@@{#1}}}
\def\frak@@#1{\noaccents@\fam\euffam#1}

\edef\msx@{\hexnumber@\msxfam}
\edef\msy@{\hexnumber@\msyfam}

\mathchardef\boxdot="2\msx@00
\mathchardef\boxplus="2\msx@01
\mathchardef\boxtimes="2\msx@02
\mathchardef\square="0\msx@03
\mathchardef\blacksquare="0\msx@04
\mathchardef\centerdot="2\msx@05
\mathchardef\lozenge="0\msx@06
\mathchardef\blacklozenge="0\msx@07
\mathchardef\circlearrowright="3\msx@08
\mathchardef\circlearrowleft="3\msx@09
\mathchardef\rightleftharpoons="3\msx@0A
\mathchardef\leftrightharpoons="3\msx@0B
\mathchardef\boxminus="2\msx@0C
\mathchardef\Vdash="3\msx@0D
\mathchardef\Vvdash="3\msx@0E
\mathchardef\vDash="3\msx@0F
\mathchardef\twoheadrightarrow="3\msx@10
\mathchardef\twoheadleftarrow="3\msx@11
\mathchardef\leftleftarrows="3\msx@12
\mathchardef\rightrightarrows="3\msx@13
\mathchardef\upuparrows="3\msx@14
\mathchardef\downdownarrows="3\msx@15
\mathchardef\upharpoonright="3\msx@16

\mathchardef\downharpoonright="3\msx@17
\mathchardef\upharpoonleft="3\msx@18
\mathchardef\downharpoonleft="3\msx@19
\mathchardef\rightarrowtail="3\msx@1A
\mathchardef\leftarrowtail="3\msx@1B
\mathchardef\leftrightarrows="3\msx@1C
\mathchardef\rightleftarrows="3\msx@1D
\mathchardef\Lsh="3\msx@1E
\mathchardef\Rsh="3\msx@1F
\mathchardef\rightsquigarrow="3\msx@20
\mathchardef\leftrightsquigarrow="3\msx@21
\mathchardef\looparrowleft="3\msx@22
\mathchardef\looparrowright="3\msx@23
\mathchardef\circeq="3\msx@24
\mathchardef\succsim="3\msx@25
\mathchardef\gtrsim="3\msx@26
\mathchardef\gtrapprox="3\msx@27
\mathchardef\multimap="3\msx@28
\mathchardef\therefore="3\msx@29
\mathchardef\because="3\msx@2A
\mathchardef\doteqdot="3\msx@2B

\mathchardef\triangleq="3\msx@2C
\mathchardef\precsim="3\msx@2D
\mathchardef\lesssim="3\msx@2E
\mathchardef\lessapprox="3\msx@2F
\mathchardef\eqslantless="3\msx@30
\mathchardef\eqslantgtr="3\msx@31
\mathchardef\curlyeqprec="3\msx@32
\mathchardef\curlyeqsucc="3\msx@33
\mathchardef\preccurlyeq="3\msx@34
\mathchardef\leqq="3\msx@35
\mathchardef\leqslant="3\msx@36
\mathchardef\lessgtr="3\msx@37
\mathchardef\backprime="0\msx@38
\mathchardef\risingdotseq="3\msx@3A
\mathchardef\fallingdotseq="3\msx@3B
\mathchardef\succcurlyeq="3\msx@3C
\mathchardef\geqq="3\msx@3D
\mathchardef\geqslant="3\msx@3E
\mathchardef\gtrless="3\msx@3F
\mathchardef\sqsubset="3\msx@40
\mathchardef\sqsupset="3\msx@41
\mathchardef\vartriangleright="3\msx@42
\mathchardef\vartriangleleft="3\msx@43
\mathchardef\trianglerighteq="3\msx@44
\mathchardef\trianglelefteq="3\msx@45
\mathchardef\bigstar="0\msx@46
\mathchardef\between="3\msx@47
\mathchardef\blacktriangledown="0\msx@48
\mathchardef\blacktriangleright="3\msx@49
\mathchardef\blacktriangleleft="3\msx@4A
\mathchardef\vartriangle="0\msx@4D
\mathchardef\blacktriangle="0\msx@4E
\mathchardef\triangledown="0\msx@4F
\mathchardef\eqcirc="3\msx@50
\mathchardef\lesseqgtr="3\msx@51
\mathchardef\gtreqless="3\msx@52
\mathchardef\lesseqqgtr="3\msx@53
\mathchardef\gtreqqless="3\msx@54
\mathchardef\Rrightarrow="3\msx@56
\mathchardef\Lleftarrow="3\msx@57
\mathchardef\veebar="2\msx@59
\mathchardef\barwedge="2\msx@5A
\mathchardef\doublebarwedge="2\msx@5B
\mathchardef\angle="0\msx@5C
\mathchardef\measuredangle="0\msx@5D
\mathchardef\sphericalangle="0\msx@5E
\mathchardef\varpropto="3\msx@5F
\mathchardef\smallsmile="3\msx@60
\mathchardef\smallfrown="3\msx@61
\mathchardef\Subset="3\msx@62
\mathchardef\Supset="3\msx@63
\mathchardef\Cup="2\msx@64

\mathchardef\Cap="2\msx@65

\mathchardef\curlywedge="2\msx@66
\mathchardef\curlyvee="2\msx@67
\mathchardef\leftthreetimes="2\msx@68
\mathchardef\rightthreetimes="2\msx@69
\mathchardef\subseteqq="3\msx@6A
\mathchardef\supseteqq="3\msx@6B
\mathchardef\bumpeq="3\msx@6C
\mathchardef\Bumpeq="3\msx@6D
\mathchardef\lll="3\msx@6E

\mathchardef\ggg="3\msx@6F

\mathchardef\circledS="0\msx@73
\mathchardef\pitchfork="3\msx@74
\mathchardef\dotplus="2\msx@75
\mathchardef\backsim="3\msx@76
\mathchardef\backsimeq="3\msx@77
\mathchardef\complement="0\msx@7B
\mathchardef\intercal="2\msx@7C
\mathchardef\circledcirc="2\msx@7D
\mathchardef\circledast="2\msx@7E
\mathchardef\circleddash="2\msx@7F
\def\ulcorner{\delimiter"4\msx@70\msx@70 }
\def\urcorner{\delimiter"5\msx@71\msx@71 }
\def\llcorner{\delimiter"4\msx@78\msx@78 }
\def\lrcorner{\delimiter"5\msx@79\msx@79 }
\def\yen{\mathhexbox\msx@55 }
\def\checkmark{\mathhexbox\msx@58 }
\def\circledR{\mathhexbox\msx@72 }
\def\maltese{\mathhexbox\msx@7A }
\mathchardef\lvertneqq="3\msy@00
\mathchardef\gvertneqq="3\msy@01
\mathchardef\nleq="3\msy@02
\mathchardef\ngeq="3\msy@03
\mathchardef\nless="3\msy@04
\mathchardef\ngtr="3\msy@05
\mathchardef\nprec="3\msy@06
\mathchardef\nsucc="3\msy@07
\mathchardef\lneqq="3\msy@08
\mathchardef\gneqq="3\msy@09
\mathchardef\nleqslant="3\msy@0A
\mathchardef\ngeqslant="3\msy@0B
\mathchardef\lneq="3\msy@0C
\mathchardef\gneq="3\msy@0D
\mathchardef\npreceq="3\msy@0E
\mathchardef\nsucceq="3\msy@0F
\mathchardef\precnsim="3\msy@10
\mathchardef\succnsim="3\msy@11
\mathchardef\lnsim="3\msy@12
\mathchardef\gnsim="3\msy@13
\mathchardef\nleqq="3\msy@14
\mathchardef\ngeqq="3\msy@15
\mathchardef\precneqq="3\msy@16
\mathchardef\succneqq="3\msy@17
\mathchardef\precnapprox="3\msy@18
\mathchardef\succnapprox="3\msy@19
\mathchardef\lnapprox="3\msy@1A
\mathchardef\gnapprox="3\msy@1B
\mathchardef\nsim="3\msy@1C
\mathchardef\ncong="3\msy@1D

\mathchardef\varsubsetneq="3\msy@20
\mathchardef\varsupsetneq="3\msy@21
\mathchardef\nsubseteqq="3\msy@22
\mathchardef\nsupseteqq="3\msy@23
\mathchardef\subsetneqq="3\msy@24
\mathchardef\supsetneqq="3\msy@25
\mathchardef\varsubsetneqq="3\msy@26
\mathchardef\varsupsetneqq="3\msy@27
\mathchardef\subsetneq="3\msy@28
\mathchardef\supsetneq="3\msy@29
\mathchardef\nsubseteq="3\msy@2A
\mathchardef\nsupseteq="3\msy@2B
\mathchardef\nparallel="3\msy@2C
\mathchardef\nmid="3\msy@2D
\mathchardef\nshortmid="3\msy@2E
\mathchardef\nshortparallel="3\msy@2F
\mathchardef\nvdash="3\msy@30
\mathchardef\nVdash="3\msy@31
\mathchardef\nvDash="3\msy@32
\mathchardef\nVDash="3\msy@33
\mathchardef\ntrianglerighteq="3\msy@34
\mathchardef\ntrianglelefteq="3\msy@35
\mathchardef\ntriangleleft="3\msy@36
\mathchardef\ntriangleright="3\msy@37
\mathchardef\nleftarrow="3\msy@38
\mathchardef\nrightarrow="3\msy@39
\mathchardef\nLeftarrow="3\msy@3A
\mathchardef\nRightarrow="3\msy@3B
\mathchardef\nLeftrightarrow="3\msy@3C
\mathchardef\nleftrightarrow="3\msy@3D
\mathchardef\divideontimes="2\msy@3E
\mathchardef\varnothing="0\msy@3F
\mathchardef\nexists="0\msy@40
\mathchardef\mho="0\msy@66
\mathchardef\eth="0\msy@67
\mathchardef\eqsim="3\msy@68
\mathchardef\beth="0\msy@69
\mathchardef\gimel="0\msy@6A
\mathchardef\daleth="0\msy@6B
\mathchardef\lessdot="3\msy@6C
\mathchardef\gtrdot="3\msy@6D
\mathchardef\ltimes="2\msy@6E
\mathchardef\rtimes="2\msy@6F
\mathchardef\shortmid="3\msy@70
\mathchardef\shortparallel="3\msy@71
\mathchardef\smallsetminus="2\msy@72
\mathchardef\thicksim="3\msy@73
\mathchardef\thickapprox="3\msy@74
\mathchardef\approxeq="3\msy@75
\mathchardef\succapprox="3\msy@76
\mathchardef\precapprox="3\msy@77
\mathchardef\curvearrowleft="3\msy@78
\mathchardef\curvearrowright="3\msy@79
\mathchardef\digamma="0\msy@7A
\mathchardef\varkappa="0\msy@7B
\mathchardef\hslash="0\msy@7D
\mathchardef\hbar="0\msy@7E
\mathchardef\backepsilon="3\msy@7F
\def\Bbb{\ifmmode\let\next\Bbb@\else
 \def\next{\errmessage{Use \string\Bbb\space only in math mode}}\fi\next}
\def\Bbb@#1{{\Bbb@@{#1}}}
\def\Bbb@@#1{\fam\msyfam#1}

\catcode`\@=12